\documentclass{amsart}
\usepackage{amssymb}

\usepackage{amscd}

\input{tcilatex}
\input tcilatex

\begin{document}
\author{Marek Kara\'{s}}
\title{Tame automorphisms of $\Bbb{C}^{3}$ with multidegree of the form $%
(p_{1},p_{2},d_{3})$}
\date{}
\maketitle

\begin{abstract}
Let $d_{3}\geq p_{2}>p_{1}\geq 3$ be integers such that $p_{1},p_{2}$ are
prime numbers. In this paper we show that the sequence $(p_{1},p_{2},d_{3})$
is the multidegree of some tame automorphisms of $\Bbb{C}^{3}$ if and only
if $d_{3}\in p_{1}\Bbb{N}+p_{2}\Bbb{N},$ i.e. if and only if $d_{3}$ is a
linear combination of $p_{1}$ and $p_{2}$ with coefficients in $\Bbb{N}.$
\end{abstract}

\section{Introduction}

Let $F=(F_{1},\ldots ,F_{n}):\Bbb{C}^{n}\rightarrow \Bbb{C}^{n}$ be any
polynomial mapping. By multidegree, denoted $\limfunc{mdeg}F,$ we call the
sequence of positive integers $(\deg F_{1},\ldots ,\deg F_{n}).$ In
dimension 2 there is a complete characterization of the sequences $%
(d_{1},d_{2})$ such that there is a polynomial automorphism $F:\Bbb{C}%
^{2}\rightarrow \Bbb{C}^{2}$ with $\limfunc{mdeg}F=(d_{1},d_{2}).$ This
characterization is a consequence of the Jung van der Kullk theorem \cite
{Jung,Kulk}. Moreover in \cite{Karas} it was proven, among other things,
that there is no tame automorphism of $\Bbb{C}^{3}$ with multidegree $%
(3,4,5),(3,5,7),(4,5,7)$ and $(4,5,11).$

Recall that a tame automorphism is, by definition, a composition of linear
automorphisms and triangular automorphisms, where a triangular automorphism
is a mapping of the following form 
\begin{equation*}
T:\Bbb{C}^{n}\ni \left\{ 
\begin{array}{l}
x_{1} \\ 
x_{2} \\ 
\vdots \\ 
x_{n}
\end{array}
\right\} \mapsto \left\{ 
\begin{array}{l}
x_{1} \\ 
x_{2}+f_{2}(x_{1}) \\ 
\vdots \\ 
x_{n}+f_{n}(x_{1},\ldots ,x_{n-1})
\end{array}
\right\} \in \Bbb{C}^{n}.
\end{equation*}
By $\limfunc{Tame}(\Bbb{C}^{n})$ we will denote the group of all tame
automorphimsm of $\Bbb{C}^{n},$ and by $\limfunc{mdeg}$ the mapping from the
set of all polynomial endomorphisms of $\Bbb{C}^{n}$ into the set $\Bbb{N}%
^{n}.$ Using this notation, above mentioned facts can be written as follows $%
(3,4,5),(3,5,7),(4,5,7),(4,5,11)\notin \limfunc{mdeg}(\limfunc{Tame}(\Bbb{C}%
^{3})).$

In this paper we make a next progress in the investigation of the set $%
\limfunc{mdeg}(\limfunc{Tame}(\Bbb{C}^{3})).$ Namely we show the following
theorem.

\begin{theorem}
\label{tw_p1_p2_d3}Let $d_{3}\geq p_{2}>p_{1}\geq 3$ be positive integers.
If $p_{1}$ and $p_{2}$ are prime numbers, then $(p_{1},p_{2},d_{3})\in 
\limfunc{mdeg}(\limfunc{Tame}(\Bbb{C}^{3}))$ if and only if $d_{3}\in p_{1}%
\Bbb{N}+p_{2}\Bbb{N},$ i.e. if and only if $d_{3}$ is a linear combination
of $p_{1}$ and $p_{2}$ with coefficients in $\Bbb{N}.$
\end{theorem}

Notice that for all permutation $\sigma $ of the set $\{1,2,3\},$ $%
(d_{1},d_{2},d_{3})\in \limfunc{mdeg}(\limfunc{Tame}(\Bbb{C}^{3}))$ if and
only if $(d_{\sigma (1)},d_{\sigma (2)},d_{\sigma (3)})\in \limfunc{mdeg}(%
\limfunc{Tame}(\Bbb{C}^{3})).$ Since, also, $(d_{1},d_{2},d_{3})\in \limfunc{%
mdeg}(\limfunc{Tame}(\Bbb{C}^{3}))$ if $d_{1}=d_{2}$ (by Proposition \ref
{prop_sum_d_i} below) and $(2,d_{2},d_{3})\in \limfunc{mdeg}(\limfunc{Tame}(%
\Bbb{C}^{3}))$ for all $d_{3}\geq d_{2}\geq 2$ (\cite{Karas} Corollary 2.3),
then the assumption $d_{3}\geq p_{2}>p_{1}\geq 3$ is not restrictive.

\section{Proof of the theorem}

First, we recall one classical result (due to Sylvester) from the number
theory, particularly from so-called coin problem or Frobenius problem \cite
{Brauer}.

\begin{theorem}
\label{tw_sywester}If $a,b$ are positive integers such that $\gcd (a,b)=1,$
then for every integer $k\geq (a-1)(b-1)$ there are $k_{1},k_{2}\in \Bbb{N}$
such that 
\begin{equation*}
k=k_{1}a+k_{2}b.
\end{equation*}
Moreover $(a-1)(b-1)-1\notin a\Bbb{N}+b\Bbb{N}.$
\end{theorem}

In the proof we will, also, use the following proposition.

\begin{proposition}
\textit{(\cite{Karas}, Proposition 2.2) }\label{prop_sum_d_i}If for a
sequence of integers $1\leq d_{1}\leq \ldots \leq d_{n}$ there is $i\in
\{1,\ldots ,n\}$ such that 
\begin{equation*}
d_{i}=\sum_{j=1}^{i-1}k_{j}d_{j}\qquad \text{with }k_{j}\in \Bbb{N},
\end{equation*}
then there exists a tame automorphism $F$ of $\Bbb{C}^{n}$ with $\limfunc{%
mdeg}F=(d_{1},\ldots ,d_{n}).$
\end{proposition}

By the above proposition, in order to prove Theorem \ref{tw_p1_p2_d3}, it is
enough to show that if $d_{3}\notin p_{1}\Bbb{N}+p_{2}\Bbb{N},$ then $%
(p_{1},p_{2},d_{3})\notin \limfunc{mdeg}(\limfunc{Tame}(\Bbb{C}^{3})).$

In the proof of the above implication we will use some results and notions
from the paper of Shestakov and Umirbayev \cite{sh umb1,sh umb2}.

The first one is the following

\begin{definition}
\textit{(\cite{sh umb1}, Definition 1) }\label{def_*-red}A pair $f,g\in
k[X_{1},\ldots ,X_{n}]$ is called *-reduced if\newline
(i) $f,g$ are algebraically independent;\newline
(ii) $\overline{f},\overline{g}$ are algebraically dependent, where $%
\overline{h}$ denotes the highest homogeneous part of $h$;\newline
(iii) $\overline{f}\notin k[\overline{g}]$ and $\overline{g}\notin k[%
\overline{f}].$
\end{definition}

\begin{definition}
\textit{(\cite{sh umb1}, Definition 1) }Let $f,g\in k[X_{1},\ldots ,X_{n}]$
be a *-reduced pair with $\deg f<\deg g.$ Put $p=\frac{\deg f}{\gcd (\deg
f,\deg g)}.$ In this situation the pair $f,g$ is called $p-$reduced pair.
\end{definition}

\begin{theorem}
\textit{(\cite{sh umb1}, Theorem 2)}\label{tw_deg_g_fg} Let $f,g\in
k[X_{1},\ldots ,X_{n}]$ be a $p-$reduced pair, and let $G(x,y)\in k[x,y]$
with $\deg _{y}G(x,y)=pq+r,0\leq r<p.$ Then 
\begin{equation*}
\deg G(f,g)\geq q\left( p\deg g-\deg g-\deg f+\deg [f,g]\right) +r\deg g.
\end{equation*}
\end{theorem}

In the above theorem $[f,g]$ means the Poisson bracket of $f$ and $g,$ but
for us it is only important that 
\begin{equation*}
\deg [f,g]=2+\underset{1\leq i<j\leq n}{\max }\deg \left( \frac{\partial f}{%
\partial x_{i}}\frac{\partial g}{\partial x_{j}}-\frac{\partial f}{\partial
x_{j}}\frac{\partial g}{\partial x_{i}}\right)
\end{equation*}
if $f,g$ are algebraically independent, and $\deg [f,g]=0$ if $f,g$ are
algebraically dependent.

Notice, also, that the estimation from Theorem \ref{tw_deg_g_fg} is true
even if the condition (ii) of Definition \ref{def_*-red} is not satisfied.
Indeed, if $G(x,y)=\sum_{i,j}a_{i,j}x^{i}y^{j},$ then, by the algebraic
independence of $\overline{f}$ and $\overline{g}$ we have: 
\begin{eqnarray*}
\deg G(f,g) &=&\underset{i,j}{\max }\deg (a_{i,j}f^{i}g^{j})\geq \deg
_{y}G(x,y)\cdot \deg g= \\
&=&(qp+r)\deg g\geq q(p\deg g-\deg f-\deg g+\deg [f,g])+r\deg g.
\end{eqnarray*}
The last inequality is a consequence of the fact that $\deg [f,g]\leq \deg
f+\deg g.$

We will also use the following theorem.

\begin{theorem}
\label{tw_type_1-4}\textit{(\cite{sh umb1}, Theorem 3) }Let $%
F=(F_{1},F_{2},F_{3})\,$be a tame automorphism of $\Bbb{C}^{3}.$ If $\deg
F_{1}+\deg F_{2}+\deg F_{3}>3$ (in other words if $F$ is not a linear
automorphism), then $F$ admits either an elementary reduction or a reduction
of types I-IV (see \cite{sh umb1} Definitions 2-4).
\end{theorem}

Let us, also, recall that an automorphism $F=(F_{1},F_{2},F_{3})$ admits an
elementary reduction if there exists a polynomial $g\in \Bbb{C}[x,y]$ and a
permutation $\sigma $ of the set $\{1,2,3\}$ such that $\deg (F_{\sigma
(1)}-g(F_{\sigma (2)},F_{\sigma (3)}))<\deg F_{\sigma (1)}.$

\begin{proof}
\textit{(of Theorem \ref{tw_p1_p2_d3}) }Assume that $F=(F_{1},F_{2},F_{3})$
is an automorphism of $\Bbb{C}^{3}$ such that $\limfunc{mdeg}%
F=(p_{1},p_{2},d_{3}).$ Assume, also, that $d_{3}\notin p_{1}\Bbb{N}+p_{2}%
\Bbb{N}.$ By Theorem \ref{tw_sywester} we have: 
\begin{equation}
d_{3}<(p_{1}-1)(p_{2}-1).  \label{d_3 male}
\end{equation}
First of all we show that this hypothetical automorphism $F$ does not admit
reductions of type I-IV.

By the definitions of reductions of types I-IV (see \cite{sh umb1}
Definitions 2-4), if $F=(F_{1},F_{2},F_{3})$ admits a reduction of these
types, then $2|\deg F_{i}$ for some $i\in \{1,2,3\}.$ Thus if $d_{3}$ is
odd, then $F$ does not admit a reduction of types I-IV. Assume that $%
d_{3}=2n $ for some positive integers $n.$

If we assume that $F$ admits a reduction of type I or II, then by the
definition (see \cite{sh umb1} Definition 2 and 3) we have $p_{1}=sn$ or $%
p_{2}=sn$ for some odd $s\geq 3.$ Since $p_{1},p_{2}\leq d_{3}=2n<sn,$ then
we obtain a contradiction.

And, if we assume that $F$ admits a reduction of type III or IV, then by the
definition (see \cite{sh umb1} Definition 4) we have: 
\begin{equation*}
n<p_{1}\leq \tfrac{3}{2}n,\qquad p_{2}=3n
\end{equation*}
or 
\begin{equation*}
p_{1}=\tfrac{3}{2}n,\qquad \tfrac{5}{2}n<p_{2}\leq 3n.
\end{equation*}
Since $p_{1},p_{2}\leq d_{3}=2n<\tfrac{5}{2}n,3n,$ then we obtain a
contradiction. Thus we have proved that our hypothetical automorphism $F$
does not admit a reduction of types I-IV.

Now we will show that it, also, does not admit an elementary reduction.

Assume, by a contrary, that 
\begin{equation*}
(F_{1},F_{2},F_{3}-g(F_{1},F_{2})),
\end{equation*}
where $g\in k[x,y],$ is an elementary reduction of $(F_{1},F_{2},F_{3}).$
Then we have $\deg g(F_{1},F_{2})=\deg F_{3}=d_{3}.$ But, by Theorem \ref
{tw_deg_g_fg}, we have 
\begin{equation*}
\deg g(F_{1},F_{2})\geq q(p_{1}p_{2}-p_{1}-p_{2}+\deg [F_{1},F_{2}])+rp_{2},
\end{equation*}
where $\deg _{y}g(x,y)=qp_{1}+r$ with $0\leq r<p_{1}.\,$ Since $F_{1},F_{2}$
are algebraically independent, $\deg [F_{1},F_{2}]\geq 2$ and then 
\begin{equation*}
p_{1}p_{2}-p_{1}-p_{2}+\deg [F_{1},F_{2}]\geq
p_{1}p_{2}-p_{1}-p_{2}+2>(p_{1}-1)(p_{2}-1).
\end{equation*}
This and (\ref{d_3 male}) follows that $q=0,$ and that: 
\begin{equation*}
g(x,y)=\sum_{i=0}^{p_{1}-1}g_{i}(x)y^{i}.
\end{equation*}
Since $\func{lcm}(p_{1},p_{2})=p_{1}p_{2},$ then the sets 
\begin{equation*}
p_{1}\Bbb{N},p_{2}+p_{1}\Bbb{N},\ldots ,(p_{1}-1)p_{2}+p_{1}\Bbb{N}
\end{equation*}
are disjoint. This follows that: 
\begin{equation*}
\deg \left( \sum_{i=0}^{p_{1}-1}g_{i}(F_{1})F_{2}^{i}\right) =\underset{%
i=0,\ldots ,p_{1}-1}{\max }\left( \deg F_{1}\deg g_{i}+i\deg F_{2}\right) .
\end{equation*}
Since, also, 
\begin{equation*}
d_{3}\notin \bigcup_{r=0}^{p_{1}-1}\left( rp_{2}+p_{1}\Bbb{N}\right)
\end{equation*}
(because $d_{3}\notin p_{1}\Bbb{N}+p_{2}\Bbb{N}$), then it is easy to see
that 
\begin{equation*}
\deg \left( \sum_{i=0}^{p_{1}-1}g_{i}(F_{1})F_{2}^{i}\right) =d_{3}
\end{equation*}
is impossible.

Now, assume that 
\begin{equation*}
(F_{1},F_{2}-g(F_{1},F_{3}),F_{3}),
\end{equation*}
where $g\in k[x,y],$ is an elementary reduction of $(F_{1},F_{2},F_{3}).$
Since $d_{3}\notin p_{1}\Bbb{N}+p_{2}\Bbb{N},p_{1}\nmid d_{3}$ and $\gcd
(p_{1},d_{3})=1.$ This means, by Theorem \ref{tw_deg_g_fg}, that 
\begin{equation*}
\deg g(F_{1},F_{3})\geq q(p_{1}d_{3}-d_{3}-p_{1}+\deg [F_{1},F_{3}])+rd_{3},
\end{equation*}
where $\deg _{y}g(x,y)=qp_{1}+r$ with $0\leq r<p_{1}.$ Since $%
p_{1}d_{3}-d_{3}-p_{1}+\deg [F_{1},F_{3}]\geq p_{1}d_{3}-2d_{3}\geq
d_{3}>p_{2}$ and since we want to have $\deg g(F_{1},F_{3})=p_{2},$ then $%
q=r=0.$ This means that $g(x,y)=g(x).$ But since $p_{2}\notin p_{1}\Bbb{N},$
then equality $\deg g(F_{1})=p_{2}$ is impossible.

Finally, if we assume that $(F_{1}-g(F_{2},F_{3}),F_{2},F_{3})$ is an
elementary reduction of $(F_{1},F_{2},F_{3}),$ then in the same way as in
the previous case we obtain a contradiction.
\end{proof}

\section{Some consequences}

\begin{theorem}
Let $p_{2}>3$ be prime number and $d_{3}\geq p_{2}$ be integer. Then $%
(3,p_{2},d_{3})\in \limfunc{mdeg}(\limfunc{Tame}(\Bbb{C}^{3}))$ if and only
if $d_{3}\notin \{2p_{2}-3k\ |\ k=1,\ldots ,\left[ \frac{p_{2}}{3}\right]
\}. $
\end{theorem}

\begin{proof}
Since $p_{2}>3$ is a prime number, $p_{2}\equiv 1(\func{mod}3)$ or $%
p_{2}\equiv 2(\func{mod}3).$ Let $r\in \{1,2\}$ be such that $p_{2}\equiv r(%
\func{mod}3).$ It is easy to see that if $d_{3}\geq p_{2}$ and $d_{3}\equiv
0(\func{mod}3)$ or $d_{3}\equiv r(\func{mod}3),$ then $d_{3}\in 3\Bbb{N}%
+p_{2}\Bbb{N}.$ Thus, by Theorem \ref{tw_sywester} 
\begin{equation*}
2(p_{2}-1)-1\neq 0,r(\func{mod}3).
\end{equation*}
Take any $d_{3}$ such that $p_{2}\leq d_{3}\leq 2p_{2}-3$ and $d_{3}\neq 0,r(%
\func{mod}3).$ Since $d_{3}\leq 2p_{2}-3$ and $d_{3}\equiv 2p_{2}-3(\func{mod%
}3),$ then one can see that $d_{3}\notin 3\Bbb{N}+p_{2}\Bbb{N},$ because in
other case we would have $2p_{2}-3\in 3\Bbb{N}+p_{2}\Bbb{N}$ which is a
contradiction with Theorem \ref{tw_sywester}. Thus 
\begin{eqnarray*}
\{d_{3} &\in &\Bbb{N\ }|\ d_{3}\geq p_{2},d_{3}\notin 3\Bbb{N}+p_{2}\Bbb{N}%
\}= \\
&=&\{d_{3}\in \Bbb{N\ }|\ p_{2}\leq d_{3}\leq 2p_{2}-3,d_{3}\equiv 2p_{2}-3%
\Bbb{(}\func{mod}3\Bbb{)}\} \\
&=&\{2p_{2}-3k\ |\ k=1,\ldots ,\left[ \tfrac{p_{2}}{3}\right] \}
\end{eqnarray*}
\end{proof}

\begin{theorem}
(a) If $d_{3}\geq 7,$ then $(5,7,d_{3})\notin \limfunc{mdeg}(\limfunc{Tame}(%
\Bbb{C}^{3}))$ if and only if 
\begin{equation*}
d_{3}\neq 8,9,11,13,16,18,21,23.
\end{equation*}
\newline
(b) If $d_{3}\geq 11,$ then $(5,11,d_{3})\notin \limfunc{mdeg}(\limfunc{Tame}%
(\Bbb{C}^{3}))$ if and only if 
\begin{equation*}
d_{3}\neq 12,13,14,17,18,19,23,24,28,29,34,39.
\end{equation*}
\newline
(c) If $d_{3}\geq 13,$ then $(5,13,d_{3})\notin \limfunc{mdeg}(\limfunc{Tame}%
(\Bbb{C}^{3}))$ if and only if 
\begin{equation*}
d_{3}\neq 14,16,17,19,21,22,24,27,29,32,34,37,42,47.
\end{equation*}
\newline
(d) If $d_{3}\geq 11,$ then $(7,11,d_{3})\notin \limfunc{mdeg}(\limfunc{Tame}%
(\Bbb{C}^{3}))$ if and only if 
\begin{equation*}
d_{3}\neq 12,13,15,16,17,19,20,23,24,26,27,30,31,34,37,38,41,45,45,48,52,59.
\end{equation*}
\end{theorem}

\begin{proof}
This is a consequence of Theorem \ref{tw_sywester} and Theorem \ref
{tw_p1_p2_d3}. For example to prove (a), by Theorem \ref{tw_sywester} and
Theorem \ref{tw_p1_p2_d3} we only have to check which of the numbers $%
7,8,\ldots ,23=(5-1)(7-1)-1$ are elements of the set $5\Bbb{N}+7\Bbb{N}.$
\end{proof}

\vspace{1cm}

\textsc{Marek Kara\'{s}\newline
Instytut Matematyki\newline
Uniwersytetu Jagiello\'{n}skiego\newline
ul. \L ojasiewicza 6}\newline
\textsc{30-348 Krak\'{o}w\newline
Poland\newline
} e-mail: Marek.Karas@im.uj.edu.pl

\end{document}